\def\version{5.38}
\def\journal{{\small After submission to} TAMS}
\def\titlep{Morita isomorphism for Cuntz algebras}
\documentclass[10pt]{article}
\usepackage{graphicx,ifthen}
\usepackage{amssymb}
\usepackage{amsmath}
\usepackage{latexsym}
\usepackage{epic,eepic}

 at12pt

\newcommand{\qed}{\hbox{\rule[-2pt]{3pt}{6pt}}}
\newcommand{\qedh}{\hfill\qed \\}

\newcommand{\vv}{\vspace{.3in}}

\newtheorem{Thm}{Theorem}[section]

\newtheorem{rem}[Thm]{Remark}

\newtheorem{lem}[Thm]{Lemma}

\newtheorem{prop}[Thm]{Proposition}
\newtheorem{cor}[Thm]{Corollary}

\newcommand{\kn}{\Large\bf
$K\hspace{-.4cm} N$
\Large\bf\vv }

\def\cal#1{\mathcal #1}
\def\con{{\cal O}_{n}}
\def\coni{{\cal O}_{\infty}}

\def\pr{{\it Proof.}\quad}
\def\co#1{{\cal O}_{#1}}
\def\disp#1{{\displaystyle #1}}

\addtocounter{footnote}{1}
\def\sftt#1{
\setcounter{equation}{0}
\addtocounter{footnote}{1}
\section{#1}
}

\def\ssft#1{\subsection{#1}}

\def\cls{\quad
\clearpage
}

\begin{document}
\def\inputcls#1{\cls\input #1.txt}
\def\plan#1#2{\par\noindent\makebox[.5in][c]{#1}
\makebox[.1in][l]{$|$}
\makebox[3in][l]{#2}\\}
\newcommand{\mline}{\noindent
\thicklines
\setlength{\unitlength}{.18mm}
\begin{picture}(1000,5)
\put(0,0){\line(1,0){750}}
\end{picture}
}
\def\boxtimes{\noindent
\setlength{\unitlength}{.020918mm}
\begin{picture}(120,150)(60,0)
\thinlines
\put(0,0){
\line(1,0){100}\line(0,1){100}
}
\put(100,100){
\line(-1,0){100}\line(0,-1){100}
}
\put(40,0){$\times$}
\end{picture}
 }
\def\authorp{Katsunori  Kawamura}
\def\emailp{e-mail: kawamurakk3@gmail.com.
}
\def\addressp{{\small {\it 
College of Information Science and Engineering, 
 Ritsumeikan University,}}\\
{\small {\it 1-1-1 Noji Higashi, Kusatsu, Shiga 525-8577, Japan.}}
}
\def\ba{\mbox{\boldmath$a$}}
\def\bb{\mbox{\boldmath$b$}}
\def\bc{\mbox{\boldmath$c$}}
\def\bd{\mbox{\boldmath$d$}}
\def\be{\mbox{\boldmath$e$}}
\def\bff{\mbox{\boldmath$f$}}
\def\bl{\mbox{\boldmath$\ell$}}
\def\bm{\mbox{\boldmath$m$}}
\def\bn{\mbox{\boldmath$n$}}
\def\bp{\mbox{\boldmath$p$}}
\def\bq{\mbox{\boldmath$q$}}
\def\bu{\mbox{\boldmath$u$}}
\def\bv{\mbox{\boldmath$v$}}
\def\bw{\mbox{\boldmath$w$}}
\def\bx{\mbox{\boldmath$x$}}
\def\by{\mbox{\boldmath$y$}}
\def\bz{\mbox{\boldmath$z$}}
\def\bzero{\mbox{\boldmath$0$}}
\def\aei{almost everywhere in}
\def\ptimes{\otimes_{\varphi}}
\def\titlepage{

\noindent
{\bf 
\noindent
\thicklines
\setlength{\unitlength}{.1mm}
\begin{picture}(1000,0)(0,-100)
\put(0,0){\kn \knn\, for \journal\, Ver.\,\version}
\put(0,-50){\today,\quad {\rm file:}
 {\rm {\small \textsf{tit01.txt,\, J1.tex}}}}
\end{picture}
}
\vspace{-.5cm}
\quad\\
{\small 
\footnote{
\begin{minipage}[t]{6in}
directory: \textsf{\fileplace}, \\
file: \textsf{\incfile},\, from \startdate
\end{minipage}
}}
\quad\\
\framebox{
\begin{tabular}{ll}
\textsf{Title:} &
\begin{minipage}[t]{4in}
\titlep
\end{minipage}
\\
\textsf{Author:} &\authorp
\end{tabular}
}
{\footnotesize	
\tableofcontents }
}

\def\qmatrix#1{
\left[
\matrix{
#1
}
\right]
}
\def\dmatrix#1{
\left|
\matrix{
#1
}
\right|
}
\def\ip#1#2{\langle#1|#2\rangle}
\newcommand{\rep}{\mathop{\mathrm{Rep}}\nolimits}
\renewcommand{\hom}{\mathop{\mathrm{Hom}}\nolimits}%
\newcommand{\mor}{\mathop{\mathrm{Mor}}\nolimits}
\def\lspan#1{{\rm Lin}\langle#1\rangle}
\def\hilb{{\bf Hilb}}
%
%
%
\setcounter{section}{0}
\setcounter{footnote}{0}
\setcounter{page}{1}
\pagestyle{plain}

%
%
\title{\titlep}
\author{\authorp\thanks{\emailp}\\ \addressp}
\date{}
\maketitle
%
%
\begin{abstract}
Let ${\rm Rep}\,{\cal O}_n$ denote 
the category of all nondegenerate $^*$ representations
of the Cuntz algebra ${\cal O}_n$.
For any $2\leq n,m<\infty$,
we construct an isomorphism functor
$F_{n,m}$ from ${\rm Rep}\,{\cal O}_m$ to 
${\rm Rep}\,{\cal O}_n$ such that
\begin{enumerate}
\item
$F_{n,m}$ commutes with infinite direct sum,
\item
$F_{n,m}\circ F_{m,l}=F_{n,l}$
and $F_{m,n}=F_{n,m}^{-1}$ for any $2\leq n,m,l<\infty$,
\item
for the von Neumann algebra ${\cal N}_{\pi}$ 
generated by the image of a representation $\pi$,\
${\cal N}_{F_{n,m}(\pi)}$ and ${\cal N}_{\pi}$
are isomorphic
for any $\pi$ in ${\rm Rep}\,{\cal O}_m$, and
\item
there exists a functor $F_{\infty,n}$ from
${\rm Rep}\, {\cal O}_n$ to 
${\rm Rep}\, {\cal O}_{\infty}$
with a right inverse
such that
$F_{\infty,n}\circ F_{n,m}=F_{\infty,m}$
for any $2\leq n,m<\infty$.
\end{enumerate}
\end{abstract}

\noindent
{\bf Mathematics Subject Classifications (2020).}
16D90; 46K10; 46M15; 47A67.\\
{\bf Key words.}
Cuntz algebra, Morita equivalence, category isomorphism.

%
%
\sftt{Main theorem}
\label{section:first}
In this paper, we show a set of isomorphisms for a certain class of categories
by constructing two-sided invertible functors explicitly.

Cuntz algebras are
typical examples of 
non-type I, separable, nuclear $C^*$-algebras (\cite{C}).
From general theory (Theorem \ref{Thm:Beer})
and $K$-theory (\cite{Cuntz1981}),
it is known that
any two mutually nonisomorphic Cuntz algebras 
are Morita equivalent as $C^*$-algebras,
but not strongly Morita equivalent,
that is, they are not Morita equivalent as rings.
In this paper, we introduce a new equivalence relation of $C^*$-algebras.

Two categories ${\cal C}$ and ${\cal D}$ are said to be {\it equivalent}
({\it resp. isomorphic})
if there exists a pair $(F,G)$ of functors, 
$F:{\cal C}\to {\cal D}$ and $G:{\cal D}\to {\cal C}$ 
such that 
$F\circ G\cong id_{{\cal D}}$
and 
$G\circ F\cong id_{{\cal C}}$
({\it resp.}
$F\circ G= id_{{\cal D}}$
and 
$G\circ F= id_{{\cal C}}$, in this case, we write $F^{-1}:=G$)
where $\cong$ means isomorphism of functors 
(\cite{HS,KS,LensterE,Maclane,Mitchell}).
Hence the notion of isomorphism of categories is much stronger than 
that of equivalence of categories. 
In this sense, 
for two rings,
an isomorphism between categories of their all left modules
is stronger than their Morita equivalence.
We call this new relation {\it Morita isomorphism} of rings,
which is also an equivalence relation of rings.
For $C^*$-algebras,
``left modules" 
in the definition
is replaced with ``nondegenerate $^*$ representations".
Remark that two notions, category equivalence
and category isomorphism are often used as a same meaning
by authors.
Kiiti Morita himself used category isomorphism as the meaning of 
category equivalence (\cite[p86]{Morita1}, \cite[p452]{Morita2}).
Therefore Morita's category isomorphism
is not a category isomorphism in \cite{HS,KS,LensterE,Maclane,Mitchell}.
About examples of isomorphism and equivalence of categories,
see \cite[14.1]{HS}.

Let ${\rm Rep}\,{\cal O}_n$ denote 
the category of all nondegenerate $^*$ representations
of the Cuntz algebra ${\cal O}_n$,
which is an additive category with the direct sum as a biproduct,
the $0$ representation as a zero object,
and intertwiners as morphisms of objects
 (\cite{Maclane,Rieffel1974a}).
%
%
\begin{Thm}
\label{Thm:first}
For any $2\leq n,m<\infty$,
there exists an isomorphism functor
$F_{n,m}$ from ${\rm Rep}\,{\cal O}_m$ to 
${\rm Rep}\,{\cal O}_n$ which satisfies the following:
\begin{enumerate}
\item
$F_{n,m}$ commutes with infinite direct sum.
\item
$F_{n,m}\circ F_{m,l}=F_{n,l}$
and $F_{m,n}=F_{n,m}^{-1}$ for any $2\leq n,m,l<\infty$.
\item
For the von Neumann algebra ${\cal N}_{\pi}$ 
generated by the image of a representation $\pi$,\
${\cal N}_{F_{n,m}(\pi)}$ and ${\cal N}_{\pi}$
are isomorphic
for any $\pi$ in ${\rm Rep}\,{\cal O}_m$.
\item
There exists a functor $F_{\infty,n}$ from
${\rm Rep}\, {\cal O}_n$ to 
${\rm Rep}\, {\cal O}_{\infty}$ with a right inverse
$F_{n,\infty}$
such that
$F_{\infty,n}\circ F_{n,m}=F_{\infty,m}$
and 
$F_{n,m}\circ F_{m,\infty}=F_{n,\infty}$
for any $2\leq n,m<\infty$.
\item
For any $2\leq n,m\leq \infty$,
$\pi\in\rep\co{m}$ and $x\in\con$,
the operator
$\{F_{n,m}(\pi)\}(x)$ has a presentation as noncommutative power series
in Cuntz generators of $\pi({\cal O}_m)$ and their conjugates 
which is independent in the choice of $\pi$.
\item
Let $\hilb$ denote the category of all Hilbert spaces
{\rm (\cite{Heunen})}.
Then the forgetful functor $G_n:\rep\con\to \hilb$
satisfies
$G_n\circ F_{n,m}=G_{m}$
for any $2\leq n,m\leq \infty$,
in other words,
$F_{n,m}$ preserves any representation space.
\end{enumerate}
\end{Thm}

From Theorem \ref{Thm:first}, $\con$ and $\co{m}$
are Morita isomorphic for any $2\leq n,m<\infty$
even if they are not isomorphic. 
It is often stated that 
a categorical isomorphism seems appropriate at first glance
but it is too strong (\cite[p86]{HS}),
a categorical isomorphism
almost never appears in practice (\cite[p21]{KS}), 
it is unreasonably strict (\cite[p33]{LensterE}), and 
it is rare (\cite[p52]{Mitchell}).
In such a point of view,
results in Theorem \ref{Thm:first} are 
beyond the bounds of common sense of category theory. 

Our original aim was a study of extensions of 
representations of $\coni$ to $\con$ ($\S$ \ref{subsection:thirdfour}), but
not to find isomorphism functors
in Theorem \ref{Thm:first}.
From the study of extensions,
we found a strange similarity between
certain classes of 
representations of $\co{m}$ 
and $\con$ (see $\S$ \ref{subsection:thirdthree}).
Roughly speaking,
there exists a pair of representations 
in $\co{m}$ and $\con$ which are ``same" on $\coni$
where ``same" means not unitarily equivalent but 
exactly same on $\coni$ (see Lemma \ref{lem:transitionb}(ii)).
Theorem \ref{Thm:first} was incidentally obtained by extending 
such classes to wholes of $\rep\co{m}$'s through a trial and error process.

From Theorem \ref{Thm:first}(iii),
the von Neumann-algebraic characterization of a representation
is completely preserved by $F_{n,m}$,
that is,  the following holds.
%
%
\begin{cor}
\label{cor:first}
Fix $2\leq n,m<\infty$.
For $\pi\in\rep\co{m}$,
let  $\pi':=F_{n,m}(\pi)\in\rep\con$.
Then the following hold:
\begin{enumerate}
\item
$\pi$ is a factor representation
if and only if
$\pi'$ is a factor representation.
\item
$\pi$ is of type {\rm X}
if and only if 
$\pi'$ is of type {\rm X}
where 
the type means the type of 
representation with respect to the
Murray-von Neumann-Connes classification
and 
{\rm X} is {\rm I, II$_1$, II$_{\infty}$ or  III$_{\lambda}$}
{\rm (\cite[Chapter 5]{Dixmier}, \cite{Connes})}.
\item
$\pi$ is irreducible if and only if $\pi'$ is irreducible.
\item
$\pi=\bigoplus_{\lambda}\pi_{\lambda}$
if and only if 
$\pi'=\bigoplus_{\lambda}\pi'_{\lambda}$.
\end{enumerate}
\end{cor}

\noindent
We can omit II$_1$ in Corollary \ref{cor:first}(ii)
because there exists no representation of type II$_1$
for Cuntz algebras. 

In $\S$ \ref{section:second},
we review Cuntz algebras and their embeddings.
In $\S$ \ref{section:third},
lemmas are proved.
In $\S$ \ref{section:fourth},
we define $F_{n,m}$ and prove
Theorem \ref{Thm:first}.
In $\S$ \ref{section:fifth},
we show formulas related to Theorem \ref{Thm:first}(v).
In $\S$ \ref{section:sixth},
open problems are discussed.
In Appendix \ref{section:appone},
a short survey of Morita equivalence is given.

%
%
\sftt{An embedding of $\coni$ into $\con$}
\label{section:second}
%
%
\ssft{$\rep\con$}
\label{subsection:secondone}
For $2\leq n\leq \infty$,
let $\con$ denote the {\it Cuntz algebra}  (\cite{C}), that is,
it is a universal unital $C^*$-algebra generated by
{\it Cuntz generators} $s_1,\ldots,s_n$
which satisfy {\it Cuntz relations}
$s_i^*s_j=\delta_{ij}I$ for all $i,j=1,\ldots,n$
and $s_1s_1^*+\cdots+s_ns_n^*=I$
when $n<\infty$.
When $n=\infty$,
the last equality is replaced with the following inequality
%
%
\begin{equation}
\label{eqn:onisum}
\sum_{j=1}^{k}s_js_j^*\ \leq \ I\quad \mbox{for all }k\geq 1
\end{equation}
by definition.
Since the limit when $k\to\infty$ in 
the left hand side of 
(\ref{eqn:onisum}) does not converge with respect to
the norm topology,
the symbol $\sum_{j=1}^{\infty}s_js_j^*$ makes no sense in $\coni$.
For any $n$, $\con$ is simple, that is, there exists no
nontrivial closed two-sided ideal.
Therefore any 
nonzero $^*$ homomorphism from $\con$ to a $C^*$-algebra $A$
is injective.
Especially, any nonzero representation is injective.
Hence if operators $S_1,\ldots,S_n$ on a Hilbert space ${\cal H}$
satisfy Cuntz relations, then they give a unique 
representation of $\con$ on ${\cal H}$.

Let $\rep\con$ be as in
Theorem \ref{Thm:first}.
We write $({\cal H},\pi)\in \rep\con$
or $\pi\in \rep\con$
when 
 $({\cal H},\pi)$ is an object in $\rep\con$ for short.
Then $\rep\con$ consists of all unital representations
and the $0$ representation of $\con$.
Any nonzero object $({\cal H},\pi)\in\rep\con$
is identified with
a data $({\cal H},S_1,\ldots,S_n)$
consisting of Cuntz generators $S_1,\ldots,S_n$ on the Hilbert space ${\cal H}$.
Hence, for $\pi_1,\pi_2\in\rep\con$,
$\pi_1=\pi_2$ if and only if 
$\pi_1(s_i)=\pi_2(s_i)$ for all $i=1,\ldots,n$.
About examples of representations and states of Cuntz algebras,
see  \cite{BJ,Cuntz1980,DKS,DaPi2,DaPi3,GP08,GP15}.

%
%
\ssft{Definition}
\label{subsection:secondtwo}
Fix $2\leq n<\infty$.
Let $s_1,\ldots,s_n$ and $t_1,t_2,\ldots$ denote
Cuntz generators of $\con$ and $\coni$, respectively.
Define the embedding $f_{n,\infty}$ 
of $\coni$ into $\con$ by
%
%
\begin{equation}
\label{eqn:fminf}
f_{n,\infty}(t_{(n-1)k+i}):=s_n^k\,s_i\quad(k\geq 0,\ i=1,\ldots,n-1)
\end{equation}
where we define $s_n^0:=I$.
The embedding $f_{n,\infty}$ has been appeared in 
\cite[Definition 1.4(ii)]{GP15}
and 
\cite[Theorem 1.10]{PC01}.
Then the following hold.
%
%
\begin{lem}
\label{lem:embedding}
\begin{enumerate}
\item
For $i=1,\ldots,n-1$,
$f_{n,\infty}(t_i)=s_i$.
\item
For $j\geq 1$,
$s_nf_{n,\infty}(t_j)=f_{n,\infty}(t_{j+n-1})$.
\item
For $j\geq 1$,
%
%
\begin{equation}
\label{eqn:snf}
s_n^*f_{n,\infty}(t_{j})
=\left\{
\begin{array}{ll}
f_{n,\infty}(t_{j-n+1})\quad & (j\geq n),\\
\\
0\quad & (1\leq j\leq n-1).
\end{array}
\right.
\end{equation}
\end{enumerate}
\end{lem}
%
%
\pr
From (\ref{eqn:fminf})
and Cuntz relations of $\coni$ and $\con$,
all statements hold.
\qedh

For any $\pi\in\rep\con$,
$\pi\circ f_{n,\infty}\in\rep\coni$.
We show a property of this correspondence
as follows.
%
%
\begin{lem}
\label{lem:diff}
There exist
$\pi_1,\pi_2\in\rep\con$
such that
$\pi_1$ and $\pi_2$ are not unitarily equivalent 
but $\pi_1\circ f_{n,\infty}$
and $\pi_2\circ f_{n,\infty}$
are unitarily equivalent. 
\end{lem}
%
%
\pr
It is known that
there exist two states 
$\omega_1,\omega_2$ of $\con$
such that $\omega_1(s_n)=1$ and 
$\omega_2(s_n)=-1$.
From these equations,
both $\omega_1$ and $\omega_2$
are uniquely defined as states of $\con$.
They are called {\it Cuntz states} (\cite{Cuntz1980,GP08}).
Let $({\cal H}_i,\pi_i,\Omega_i)$
denote the Gel'fand-Naimark-Segal representation
of $\con$
defined by $\omega_i$ for $i=1,2$.
Then $\pi_1$ and $\pi_2$
are not unitarily equivalent and 
eigenequations 
$\pi_1(s_n)\Omega_1=\Omega_1$
and 
$\pi_2(s_n)\Omega_2=-\Omega_2$
hold.
Define $\eta_i:=\pi_i\circ f_{n,\infty}\in\rep\coni$.
From the definition of $f_{n,\infty}$ and the eigenequation
of $\pi_i(s_n)$,
we see $\eta_i(t_j)^*\Omega_i=0$ for all $j\geq 1$
and $i=1,2$.
Let $\rho_i:=\ip{\Omega_i}{\eta_i(\cdot)\Omega_i}$
for $i=1,2$.
Then 
we can verify that $\rho_1=\rho_2$
and 
$\Omega_i$ is also cyclic for
$({\cal H}_i,\eta_i)$ for $i=1,2$.
This implies that
$({\cal H}_1,\eta_1)$ 
and 
$({\cal H}_2,\eta_2)$ are unitarily equivalent.
\qedh

\noindent
The state $\rho_i$ of $\coni$ in the proof of 
Lemma \ref{lem:diff} has been appeared 
in \cite[EXAMPLE 5.3.27]{BR2}.

%
%
\sftt{Ingredients of functors}
\label{section:third}
In order to construct $F_{n,m}$ in Theorem \ref{Thm:first},
we construct some operators and representations.
Fix $2\leq n<\infty$.
Let $s_1,\ldots,s_n$ and $t_1,t_2,\ldots$ denote
Cuntz generators of $\con$ and $\coni$, respectively.
%
%
\ssft{Ingredients of Cuntz generators}
\label{subsection:thirdone}
For $({\cal H},\pi)\in\rep\con$ and an integer $a\geq 0$,
define operators $R_{\pi,a}$ and $Q_{\pi}$ on ${\cal H}$ by
%
%
\begin{equation}
\label{eqn:qr}
R_{\pi,a}:=
\sum_{j=1}^{\infty}\pi(f_{n,\infty}(t_{j+a}\,t_j^*)),\quad
Q_{\pi}:=R_{\pi,0}
\end{equation}
where $f_{n,\infty}$ is as in  (\ref{eqn:fminf})
and
the infinite sum is defined with respect to
the strong operator topology on ${\cal B}({\cal H})$ but not on $\pi(\con)$.
Then the following hold.
%
%
\begin{lem}
\label{lem:operators}
\begin{enumerate}
\item
\label{operatorsone}
For $a,b\geq 0$,\
$R_{\pi,a}R_{\pi,b}=R_{\pi,a+b}$.
\item
\label{operatorstwo}
For $a\geq 0$,\
$R_{\pi,a}^*R_{\pi,a}=Q_{\pi}$
and $R_{\pi,a}R_{\pi,a}^*=\sum_{j=a+1}^{\infty}
\pi(f_{m,\infty}(t_jt_j^*))$.
\item
\label{operatorsthree}
For $a\geq 0$,\
$Q_{\pi}R_{\pi,a}=R_{\pi,a}=R_{\pi,a}Q_{\pi}$.
\item
\label{operatorsfour}
For $a\geq 0$ and $j\geq 1$,\
$R_{\pi,a}\pi(f_{n,\infty}(t_{j}))=\pi(f_{n,\infty}(t_{j+a}))$
and $Q_{\pi}\pi(f_{n,\infty}(t_{j}))=\pi(f_{n,\infty}(t_{j}))$.
\item
\label{operatorsfive}
For $a\geq 0$,
$\pi(s_n)R_{\pi,a}=R_{\pi,a+n-1}$
and 
$\pi(s_n)Q_{\pi}=R_{\pi,n-1}$.
\item
\label{operatorssix}
If $\pi=\bigoplus_{\lambda}\pi_{\lambda}$, 
then $R_{\pi,a}=\sum_{\lambda}R_{\pi_{\lambda},a}$
for any $a\geq 0$.
\item
\label{operatorsseven}
For $a\geq 0$ and $j\geq 1$,
%
%
\begin{equation}
\label{eqn:pifninf}
\pi(f_{n,\infty}(t_j))^*R_{\pi,a}=
\left\{
\begin{array}{ll}
\pi(f_{n,\infty}(t_{j-a}))^*\quad &(j\geq a+1),\\
\\
0\quad & (j\leq a).
\end{array}
\right.
\end{equation}
\end{enumerate}
\end{lem}
%
%
\pr
Since $Q_{\pi}$ is a projection on ${\cal H}$
and $Q_{\pi}{\cal H}=\oplus_{j=1}^{\infty}\pi(f_{n,\infty}(t_{j})){\cal H}$,
we obtain
%
%
\begin{equation}
\label{eqn:qdeco}
{\cal H}=(I-Q_{\pi}){\cal H}\ \oplus \
\bigoplus_{j=1}^{\infty}\pi(f_{n,\infty}(t_{j})){\cal H}.
\end{equation}
On each direct sum component in (\ref{eqn:qdeco}),
all relations of operators can be proved by Lemma \ref{lem:embedding}.
\qedh

In general,
$Q_{\pi}$ is a projection but not the identity operator on ${\cal H}$.
For example,
if $({\cal H},\pi)\in \rep\con$
has a nonzero vector $\Omega\in {\cal H}$
such that $\pi(s_n)\Omega=\Omega$ (see the proof of Lemma \ref{lem:diff}),
then we can verify $Q_{\pi}\Omega=0$.
Hence $0\lneq Q_{\pi} \lneq I$.
On the other hand,
if $({\cal H},\pi)\in\rep\con$
has a cyclic vector $\Omega\in {\cal H}$
such that $\pi(s_1)\Omega=\Omega$,
then we can prove $Q_{\pi}=I$.
If $\pi=0$, then $Q_{\pi}=0$.

%
%
\begin{lem}
\label{lem:uoperator}
For $({\cal H},\pi)\in\rep\con$ with $\pi\ne 0$,
define the operator $U_{\pi}$ on ${\cal H}$ by
%
%
\begin{equation}
\label{eqn:upi}
U_{\pi}:=\pi(s_n)(I-Q_{\pi}).
\end{equation}
Then
\begin{enumerate}
\item
$U_{\pi}^*U_{\pi}=I-Q_{\pi}$.
\item
$U_{\pi}U_{\pi}^*=I-Q_{\pi}$.
\item
For $a\geq 0$,
$U_{\pi}^*R_{\pi,a}=0$. 
\item
For $a\geq 0$,
$(U_{\pi}+R_{\pi,a})^*(U_{\pi}+R_{\pi,a})=I$.
\item
For $a\geq 0$,
$(U_{\pi}+R_{\pi,a})(U_{\pi}+R_{\pi,a})^*
=I-Q_{\pi}+R_{\pi,a}R_{\pi,a}^*$.
\end{enumerate}
\end{lem}
%
%
\pr
(i) By definition, the statement holds.

\noindent
(ii)
From 
Lemma \ref{lem:operators}\ref{operatorsthree} and \ref{operatorsfive},
we obtain 
%
%
\begin{equation}
\label{eqn:upiupi}
\begin{array}{rl}
U_{\pi}U_{\pi}^*
=&
\pi(s_ns_n^*)-R_{\pi,n-1}R_{\pi,n-1}^*\\
\\
=&\disp{I-\sum_{i=1}^{n-1}\pi(s_is_i^*)-
\sum_{j=n}^{\infty}\pi(f_{n,\infty}(t_jt_j^*))}\\
\\
=&\disp{I-\sum_{i=1}^{n-1}\pi(f_{n,\infty}(t_it_i^*))-
\sum_{j=n}^{\infty}\pi(f_{n,\infty}(t_jt_j^*))
}\\
\\
=&I-Q_{\pi}.
\end{array}
\end{equation}

\noindent
(iii)
From 
Lemma \ref{lem:operators}\ref{operatorsthree} and \ref{operatorsfour},
the statement holds.

\noindent
(iv)
From (iii) and (i), the statement holds.

\noindent
(v)
From (iii) and (ii),
the statement holds.
\qedh

%
%
\ssft{Construction of a new representation from a given representation}
\label{subsection:thirdtwo}
Fix $2\leq n,m<\infty$.
Let $s_1,\ldots,s_n,r_1,\ldots,r_m$ and 
$t_1,t_2,\ldots$ denote
Cuntz generators of $\con,\co{m}$ and $\coni$,
respectively.
Remark that we use symbols $R_{\pi,a}$ and $Q_{\pi}$
 in (\ref{eqn:qr})
for $\pi$ in $\rep\co{m}$ but not for $\pi$ in $\rep\con$.
For $({\cal H},\pi)\in\rep\co{m}$,
define operators $\pi'(s_1),\ldots,\pi'(s_n)$ on ${\cal H}$ by
%
%
\begin{equation}
\label{eqn:pidash}
\pi'(s_i):=
\left\{
\begin{array}{ll}
\pi(f_{m,\infty}(t_i))\quad &(i=1,\ldots,n-1),\\
\\
\pi(r_m)(I-Q_{\pi})+R_{\pi,n-1}\quad& (i=n)
\end{array}
\right.
\end{equation}
where $f_{m,\infty}$ is as in (\ref{eqn:fminf}).
Remark positions of $n$ and $m$ in the definition of $\pi'(s_n)$.
For readability,
we write $\pi'(s_n)$ explicitly as follows:
%
%
\begin{equation}
\label{eqn:pidashntwo}
\pi'(s_n)=\pi(r_m)\left\{I-\sum_{j=1}^{\infty}
\pi(f_{m,\infty}(t_jt_j^*))
\right\}+\sum_{j=1}^{\infty}\pi(f_{m,\infty}(t_{j+n-1}\,t_j^*)).
\end{equation}
About the explicit description of 
$\pi'(s_i)$'s by $\pi(r_j)$'s, see $\S$ \ref{section:fifth}. 
For $\pi\in \rep\co{m}$ with $\pi\ne 0$,
we will prove that 
$\pi'(s_1),\ldots,\pi'(s_n)$ satisfy Cuntz relations of $\con$,
that is,
$\pi'$ is a representation of $\con$ on ${\cal H}$
in Lemma \ref{lem:rep}.

%
%
\begin{lem}
\label{lem:transition}
Let $\pi\in\rep\co{m}$
with $\pi'(s_i)$'s  in {\rm (\ref{eqn:pidash})}.
\begin{enumerate}
\item
For $j\geq 1$,
$\pi'(s_n)\pi(f_{m,\infty}(t_j))=
\pi(f_{m,\infty}(t_{j+n-1}))$.
\item
If $\pi=\bigoplus_{\lambda}\pi_{\lambda}$,
then $\pi'(s_i)=\sum_{\lambda}\pi'_{\lambda}(s_i)$
for all $i=1,\ldots,n$.
\end{enumerate}
\end{lem}
%
%
\pr
It suffices to consider only the case of $\pi\ne 0$.
Let $Q_{\pi}^{\perp}:=I-Q_{\pi}$.

\noindent
(i)
From Lemma \ref{lem:operators}\ref{operatorsfour},
we see
$Q_{\pi}^{\perp}\pi(f_{m,\infty}(t_j))=0$
and 
$R_{\pi,n-1}\pi(f_{m,\infty}(t_j))=\pi(f_{m,\infty}(t_{j+n-1}))$.
Hence
we obtain
%
%
\begin{equation}
\label{eqn:pidashpifm}
\begin{array}{rl}
\pi'(s_n)\pi(f_{m,\infty}(t_j))
=&\{\pi(r_m)Q_{\pi}^{\perp}+R_{\pi,n-1}\}\pi(f_{m,\infty}(t_j))
\\
=&\pi(r_m)Q_{\pi}^{\perp}\pi(f_{m,\infty}(t_j))
+R_{\pi,n-1}\pi(f_{m,\infty}(t_j))
\\
=&\pi(f_{m,\infty}(t_{j+n-1})).
\end{array}
\end{equation}

\noindent
(ii)
From Lemma \ref{lem:operators}\ref{operatorssix}, the statement holds.
\qedh

%
%
\begin{lem}
\label{lem:rep}
Let $\pi\in\rep\co{m}$.
If $\pi\ne 0$, 
then operators
$\pi'(s_1),\ldots,\pi'(s_n)$ in 
{\rm (\ref{eqn:pidash})} satisfy
Cuntz relations of $\con$, that is,
the following hold:
\begin{enumerate}
\item
For each $i,j=1,\ldots,n$,\ 
$\pi'(s_i)^*\pi'(s_j)=\delta_{ij}I$.
\item
$\pi'(s_1)\pi'(s_1)^*+\cdots+\pi'(s_n)\pi'(s_n)^*=I$.
\end{enumerate}
\end{lem}
%
%
\pr
Let $Q_{\pi}^{\perp}:=I-Q_{\pi}$.

\noindent
(i)
By definition,
we see
$\pi'(s_i)^*\pi'(s_j)=\delta_{ij}I$ for all $i,j=1,\ldots,n-1$.
For $i=1,\ldots,n-1$,
%
%
\begin{equation}
\label{eqn:pidashstar}
\begin{array}{rl}
\pi'(s_i)^*\pi'(s_n)
=&\pi(f_{m,\infty}(t_{i}))^*
\{\pi(r_m)Q_{\pi}^{\perp}+R_{\pi,n-1}\}\\
=&\pi(f_{m,\infty}(t_{i}))^*\pi(r_m)Q_{\pi}^{\perp}+
\pi(f_{m,\infty}(t_{i}))^*R_{\pi,n-1}\\
=&0\quad(\mbox{from Lemma \ref{lem:operators}(iv) and (vii)}).
\end{array}
\end{equation}
On the other hand,
$\pi'(s_n)^*\pi'(s_n)= I$
from Lemma \ref{lem:uoperator}(iv).
Hence the relations are proved.

\noindent
(ii)
By definition,
%
%
\begin{equation}
\label{eqn:sumcond}
\disp{\sum_{i=1}^n\pi'(s_i)\pi'(s_i)^*}
=\disp{
\sum_{i=1}^{n-1}\pi(f_{m,\infty}(t_i))\pi(f_{m,\infty}(t_i))^*}
+\pi'(s_n)\pi'(s_n)^*,
\end{equation}
and 
%
%
\begin{equation}
\label{eqn:upipi}
\begin{array}{rl}
\pi'(s_n)\pi'(s_n)^*
=&
(U_{\pi}+R_{\pi,n-1})
(U_{\pi}+R_{\pi,n-1})^*\\
=&I-Q_{\pi}+R_{\pi,n-1}R_{\pi,n-1}^*
\quad(\mbox{from Lemma \ref{lem:uoperator}(v)})
\\
=&
\disp{
I-Q_{\pi}+\sum_{j=n}^{\infty}\pi(f_{m,\infty}(t_jt_j^*))
\quad (\mbox{from Lemma \ref{lem:operators}(ii)})}.
\end{array}
\end{equation}
From this,
(\ref{eqn:sumcond}) and the definition of $Q_{\pi}$, 
we obtain the relation.
\qedh

\noindent
From Lemma \ref{lem:rep},
$\pi'$ in {\rm (\ref{eqn:pidash})}
is well defined in $\rep\con$.

We summarize the points so far as follows:
For any $({\cal H},\pi)\in\rep\co{m}$,
we can construct 
$({\cal H}',\pi')\in\rep\con$
for ${\cal H}':={\cal H}$
for any $2\leq n,m<\infty$,
without any additional assumption,
by constructing Cuntz generators of $\con$ on ${\cal H}$
from those of $\co{m}$ explicitly:
%
%
\begin{equation}
\label{eqn:mapfunc}
\rep\co{m}\ni ({\cal H},\pi)\mapsto 
({\cal H}',\pi')\in\rep\con,\quad
\pi'(\con)\subset \pi(\co{m})''.
\end{equation}
This result itself is meaningful as a method to
construct a new representation of $A$
from a given representation of $B$
for two distinct $C^*$-algebras $A$ and $B$
which makes a sense even if
 $\hom(A,B)=\emptyset$
and 
 $\hom(B,A)=\emptyset$
(see Proposition \ref{prop:ktheory}).
Furthermore we can prove that 
the correspondence in (\ref{eqn:mapfunc}) is invertible
for any $2\leq n,m<\infty$ (see the second equality in 
Theorem \ref{Thm:first}(ii)).
The construction $\pi\mapsto \pi'$
is not similar to any known construction 
of representation,
that is, restriction and induction (\cite{Rieffel1974b}).
The author himself has not found a reason why $\pi'$ in (\ref{eqn:pidash})
gives a representation
except the proof of Lemma \ref{lem:rep}.
Operators $R_{\pi,a}$ and $Q_{\pi}$
in (\ref{eqn:qr})
were introduced after the discovery of the formula in 
(\ref{eqn:pidashntwo}).

%
%
\ssft{Invariant properties of $\pi\mapsto \pi'$}
\label{subsection:thirdthree}
A morphism between 
two objects $({\cal H}_1,\pi_1),({\cal H}_2,\pi_2)$ in $\rep \con$
is their intertwiner (\cite{Rieffel1974a}).
Hence
the set of morphisms is given as follows:
%
%
\begin{equation}
\label{eqn:mordef}
\mor(\pi_1,\pi_2)=\{T\in {\cal B}({\cal H}_1,{\cal H}_2):
T\pi_1(x)=\pi_2(x)T\mbox{ for all }x\in \con\}
\end{equation}
where ${\cal B}({\cal H}_1,{\cal H}_2)$ denotes
the set of all bounded linear operators from 
${\cal H}_1$ to ${\cal H}_2$.
%
%
\begin{lem}
\label{lem:morlem}
Let $\pi'$ be as in {\rm (\ref{eqn:pidash})}.
For any $\pi_1,\pi_2\in\rep\co{m}$,
$\mor(\pi_1,\pi_2)\subset 
\mor(\pi_1',\pi_2')$.
\end{lem}
%
%
\pr
It suffices to show the case of $\pi\ne 0$.
For $T\in \mor(\pi_1,\pi_2)$,
$T\pi_1(x)=\pi_2(x)T$
for any $x\in \co{m}$.
Since $f_{m,\infty}(t_{j})$,
$f_{m,\infty}(t_{j})^*\in\co{m}$,
we see $T\pi_1'(s_i)=\pi_2'(s_i)T$
and $T\pi_1'(s_i^*)=\pi_2'(s_i^*)T$
for $i=1,\ldots,n-1$.
Let $I_i$ denote the identity operator on
the representation space of $\pi_i$ for $i=1,2$.
Then $TI_1=T\pi_1(I)=\pi_2(I)T=I_2T$.
From this,
%
%
\begin{equation}
\label{eqn:larget}
\begin{array}{rl}
T\pi'_1(s_n)
=&T\{\pi_1(r_m)(I_1-Q_{\pi_1})+R_{\pi_1,n-1}\}\\
=&\{\pi_2(r_m)(I_2-Q_{\pi_2})+R_{\pi_2,n-1}\}T\\
=&\pi'_2(s_n)T.
\end{array}
\end{equation}
As the same token,
we obtain
$T\pi'_1(s_n^*)=\pi'_2(s_n^*)T$.
From these,
$T\pi_1'(x)=\pi_2'(x)T$
for any $x\in \con$.
Hence $T\in \mor(\pi_1',\pi_2')$.
\qedh

\noindent
Furthermore, we can prove 
$\mor(\pi_1,\pi_2)=
\mor(\pi_1',\pi_2')$
(see the proof of Theorem \ref{Thm:first}(iii) in $\S$ \ref{section:fourth}).
%
%
\begin{lem}
\label{lem:transitionb}
Let $\pi'$ be as in 
{\rm (\ref{eqn:pidash})}.
\begin{enumerate}
\item
For any $\pi\in\rep\co{m}$,
$\pi'\circ f_{n,\infty}=\pi\circ f_{m,\infty}$.
\item
For any $\pi\in\rep\co{m}$ and $a\geq 0$,
$R_{\pi',a}=R_{\pi,a}$ and $Q_{\pi'}=Q_{\pi}$.
\item
If $n=m$, then $\pi'=\pi$.
\end{enumerate}
\end{lem}
%
%
\pr
(i)
For any $k\geq 0$ and $i=1,\ldots,n-1$,
%
%
\begin{equation}
\label{eqn:pidashinf}
\begin{array}{rl}
\pi'(f_{n,\infty}(t_{(n-1)k+i}))
=&
\pi'(s_n^ks_i)\\
=&
\pi'(s_n)^k\pi'(s_i)\\
=&
\pi'(s_n)^k\pi(f_{m,\infty}(t_i))\\
=&
\pi(f_{m,\infty}(t_{i+(n-1)k}))
\quad(\mbox{from Lemma \ref{lem:transition}(i)}).
\end{array}
\end{equation}
Hence the statement holds.

\noindent
(ii)
From (i), the first statement holds.
From this and $Q_{\pi}=R_{\pi,0}$,
the second holds.

\noindent
(iii)
For $i=1,\ldots,n-1$,
$\pi'(s_i)=\pi(f_{n,\infty}(t_i))=\pi(s_i)$
from Lemma \ref{lem:embedding}(i).
On the other hand,
we can prove $\pi'(s_n)=\pi(s_n)$
from Lemma \ref{lem:operators}\ref{operatorsfive}.
Hence the statement holds.
\qedh

\noindent
In consequence,
$R_{\pi,a}$, $Q_{\pi}$
and $\mor(\pi_1,\pi_2)$
are invariant under $\pi\mapsto \pi'$.

%
%
\ssft{Unmagnifying extensions of representations of $\coni$
to $\con$}
\label{subsection:thirdfour}
By using $f_{n,\infty}$ in 
(\ref{eqn:fminf}),
$\coni$ can be regarded as a subalgebra of $\con$.
For  $({\cal H},\pi)\in \rep\coni$,
if $({\cal K},\Pi)\in \rep\con$ satisfies
${\cal K}={\cal H}$ and 
$\Pi\circ f_{n,\infty}=\pi$,
then we call $({\cal K},\Pi)$ an {\it unmagnifying extension} of $\pi$
with respect to $f_{n,\infty}$.
We construct such an extension as follows.
For $({\cal H},\pi)\in\rep\coni$,
define $({\cal H}',\pi')\in\rep\con$ as
${\cal H}':={\cal H}$ and 
%
%
\begin{equation}
\label{eqn:dashinf}
\pi'(s_i)=
\left\{
\begin{array}{ll}
\pi(t_i)\quad&(i=1,\ldots,n-1),\\
\\
\disp{\pi(I)-\sum_{j=1}^{\infty}\pi(t_{j}t_j^*)
+\sum_{j=1}^{\infty}\pi(t_{j+n-1}\,t_j^*)}&(i=n).
\end{array}
\right.
\end{equation}
We can verify that
$\pi'(s_1),\ldots,\pi'(s_n)$ satisfy Cuntz relations of $\con$.
Remark that the definition of $\pi'(s_n)$ in
(\ref{eqn:dashinf}) is different from
that in (\ref{eqn:pidash}).
By definition, if $\pi=0$, then $\pi'=0$.
Furthermore the following holds.
%
%
\begin{lem}
\label{lem:extinf}
\begin{enumerate}
\item
For any $\pi\in\rep\coni$, $\pi'\circ f_{n,\infty}=\pi$.
\item
For any $\pi_1,\pi_2\in\rep\coni$, 
$\mor(\pi_1,\pi_2)\subset \mor(\pi_1',\pi_2')$.
\end{enumerate}
\end{lem}
%
%
\pr
(i)
It is sufficient to show the case of $\pi\ne 0$.
Then $\pi(I)=I$.
By definition,
%
%
\begin{equation}
\label{eqn:shiftinf}
\begin{array}{rl}
\pi'(s_n)\pi(t_j)
=&\disp{\left\{I-\sum_{a=1}^{\infty}\pi(t_at_a^*)+
\sum_{a=1}^{\infty}\pi(t_{a+n-1}\,t_a^*)\right\}\pi(t_j)}\\
\\
=& \pi(t_{j+n-1})
\end{array}
\end{equation}
for any $j\geq 1$.
For any $l\geq 0$ and $i=1,\ldots,n-1$,
%
%
\begin{equation}
\label{eqn:pidashcirtnm}
\begin{array}{rl}
\{\pi'\circ f_{n,\infty}\}(t_{(n-1)l+i})
=& \pi'(s_n^ls_i)\\
=& \pi'(s_n)^l\pi'(s_i)\\
=& \pi'(s_n)^l\pi(t_i)\\
=& \pi(t_{i+(n-1)l})\quad(\mbox{from (\ref{eqn:shiftinf})}).
\end{array}
\end{equation}
Hence the statement holds.

\noindent
(ii)
Along with the proof of 
Lemma \ref{lem:morlem}, the statement can be proved.
\qedh

\noindent
From Lemma \ref{lem:extinf}(i),
$\pi'$ in (\ref{eqn:dashinf})
is an unmagnifying extension of $\pi\in\rep\coni$
with respect to $f_{n,\infty}$.
In other words,
any representation of $\coni$
has an unmagnifying extension with respect to $f_{n,\infty}$.
A prototype of $\pi'$ in (\ref{eqn:dashinf}) was found before
that in (\ref{eqn:pidash}).

%
%
\sftt{Proof of Theorem \ref{Thm:first}}
\label{section:fourth}
For $\pi,\pi'$  in (\ref{eqn:pidash}),
we define the symbol $F_{n,m}(\pi)$ as 
%
%
\begin{equation}
\label{eqn:fnmdef}
F_{n,m}(\pi):=\pi'
\end{equation}
and define the function $F_{n,m}:
\mor(\pi_1,\pi_2)\to\mor(F_{n,m}(\pi_1),F_{n,m}(\pi_2))$
as $F_{n,m}(T):=T$ on morphisms.
From Lemma \ref{lem:morlem}, this is well defined.
From these, we obtain the functor
$F_{n,m}:\rep\co{m}\to \rep\con$.
We use the dash notation 
in (\ref{eqn:fnmdef})
for convenience.
\\

\noindent
{\it Proof of} Theorem \ref{Thm:first}.
(i) 
From Lemma \ref{lem:transition}(ii), 
the statement holds.

\noindent
(ii)
Let 
$s_1,\ldots,s_n$,
$r_1,\ldots,r_m$,
$u_1,\ldots,u_l$ 
and $t_1,t_2,\ldots$
denote Cuntz generators of $\con$, $\co{m}$, $\co{l}$ and $\coni$,
respectively.
For $\pi\in\rep\co{l}$, let $\eta:=F_{m,l}(\pi)$.
Then $(F_{n,m}\circ F_{m,l})(\pi)=\eta'$.
For $i=1,\ldots,n-1$,
%
%
\begin{equation}
\label{eqn:fnmcirc}
\begin{array}{rl}
\{(F_{n,m}\circ F_{m,l})(\pi)\}(s_i)
=&\eta'(s_i)\\
=&\eta(f_{m,\infty}(t_i))\\
=&\pi'(f_{m,\infty}(t_i))\\
=&\pi(f_{l,\infty}(t_i))\quad
(\mbox{from Lemma \ref{lem:transitionb}(i)})\\
=&F_{n,l}(\pi)(s_i).
\end{array}
\end{equation}
On the other hand,
%
%
\begin{equation}
\label{eqn:fnmcircb}
\begin{array}{rl}
\{(F_{n,m}\circ F_{m,l})(\pi)\}(s_n)
=&
\eta'(s_n)
=
\eta(r_m)(I-Q_{\eta})+R_{\eta,n-1}.
\end{array}
\end{equation}
Then
%
%
\begin{equation}
\label{eqn:etarm}
\begin{array}{rl}
\eta(r_m)(I-Q_{\eta})
=&
\pi'(r_m)(I-Q_{\pi'})\\
=&
\pi'(r_m)(I-Q_{\pi})\quad 
(\mbox{from Lemma \ref{lem:transitionb}(ii)})
\\
=&
\{\pi(u_l)(I-Q_{\pi})+R_{\pi,m-1}\}(I-Q_{\pi})
\\
=&
\pi(u_l)(I-Q_{\pi})+R_{\pi,m-1}(I-Q_{\pi})\\
=&
\pi(u_l)(I-Q_{\pi})\quad(\mbox{from Lemma \ref{lem:operators}(iii)}).
\end{array}
\end{equation}
From Lemma \ref{lem:transitionb}(ii),
$R_{\eta,n-1}
=
R_{\pi',n-1}
=
R_{\pi,n-1}$.
From this and (\ref{eqn:etarm}),
%
%
\begin{equation}
\label{eqn:piuli}
\mbox{(\ref{eqn:fnmcircb})}=
\pi(u_l)(I-Q_{\pi})
+R_{\pi,n-1}=F_{n,l}(\pi)(s_n).
\end{equation}
From this and (\ref{eqn:fnmcirc}),
we obtain 
$\{(F_{n,m}\circ F_{m,l})(\pi)\}(s_i)=F_{n,l}(\pi)(s_i)$ for all $i$.
Hence the first equality of functors holds.
From Lemma \ref{lem:transitionb}(iii),
$F_{n,n}=id_{\rep\con}$.
From this and the first equality,
we obtain 
$F_{n,m}\circ F_{m,n}=id_{\rep\con}$
and 
$F_{m,n}\circ F_{n,m}=id_{\rep\co{m}}$.
Hence the second equality is proved.

\noindent
(iii)
From (ii),
$\pi=(F_{m,n}\circ F_{n,m})(\pi)$
for all $\pi\in \rep\co{m}$.
Let $\pi_i\in\rep\co{m}$
and $\eta_i:=F_{n,m}(\pi_i)$ for $i=1,2$.
Then $\eta_i'=F_{m,n}(F_{n,m}(\pi_i))=\pi_i$ for $i=1,2$.
From this and Lemma \ref{lem:morlem},
%
%
\begin{equation}
\label{eqn:moreqtwo}
\mor(\pi_1,\pi_2)\subset 
\mor(\pi_1',\pi_2')
=
\mor(\eta_1,\eta_2)
\subset 
\mor(\eta_1',\eta_2')=
\mor(\pi_1,\pi_2).
\end{equation}
Hence we obtain
$\mor(\pi_1',\pi_2')=\mor(\pi_1,\pi_2)$
as a subspace of ${\cal B}({\cal H}_1,{\cal H}_2)$
where ${\cal H}_i$ denotes
the representation space of $\pi_i$ for $i=1,2$.
From this,
%
%
\begin{equation}
\label{eqn:calnfnm}
\begin{array}{rl}
{\cal N}_{F_{n,m}(\pi)}'
=&
\{F_{n,m}(\pi)(\con)\}'\\
=&\mor(F_{n,m}(\pi),
F_{n,m}(\pi))\\
=&
\mor(\pi',\pi')\\
=&
\mor(\pi,\pi)\\
=&
\pi(\co{m})'\\
=&
{\cal N}_{\pi}'.
\end{array}
\end{equation}
This implies
${\cal N}_{F_{n,m}(\pi)}={\cal N}_{\pi}$.
Hence the statement holds.

\noindent
(iv)
For $f_{n,\infty}$ in  (\ref{eqn:fminf}),
define the functor
 $F_{\infty,n}:\rep\con\to \rep\coni$ by
%
%
\begin{equation}
\label{eqn:finfncal}
F_{\infty,n}({\cal H},\pi):=({\cal H},\ \pi\circ f_{n,\infty}).
\end{equation}
Then we see $\mor(\pi_1,\pi_2)\subset 
\mor(F_{\infty,n}(\pi_1),F_{\infty,n}(\pi_2))$,
and define $F_{\infty,n}(T):=T$ on morphisms.
For $\pi\in\rep\co{m}$,
$k\geq 0$ and $i=1,\ldots,n-1$,
we obtain
%
%
\begin{equation}
\label{eqn:finfcrcb}
\begin{array}{rl}
\{(F_{\infty,n}\circ F_{n,m})(\pi)\}(t_{(n-1)k+i})
=&
F_{n,m}(\pi)(f_{n,\infty}(t_{(n-1)k+i}))\\
=&
F_{n,m}(\pi)(s_n^ks_i)\\
=&
\pi'(s_n^ks_i)\\
=&
\pi'(s_n)^k\pi'(s_i)\\
=&
\pi'(s_n)^k\pi(f_{m,\infty}(t_i))\\
=&
\pi(f_{m,\infty}(t_{i+(n-1)k}))\quad 
(\mbox{from Lemma \ref{lem:transition}(i)})\\
=&F_{\infty,m}(\pi)(t_{(n-1)k+i}).
\end{array}
\end{equation}
Hence the equality of functors is proved.

Next we construct a right inverse of $F_{\infty,n}$.
For $\pi'$ in (\ref{eqn:dashinf}),
%
%
\begin{equation}
\label{eqn:pidashfour}
F_{n,\infty}(\pi):=\pi'
\end{equation}
defines a functor from $\rep\coni$ to $\rep\con$
such that $F_{n,\infty}(T)=T$ on morphisms.
From Lemma \ref{lem:extinf}, $F_{n,\infty}$ is well defined and
$F_{\infty,n}\circ F_{n,\infty}=id_{\rep\coni}$.
Along with (\ref{eqn:finfcrcb}),
the equality for $F_{\infty,n}$'s can be verified.

\noindent
(v)
It suffices to show the case that $x$ is a Cuntz generator.
From (\ref{eqn:fminf}) and (\ref{eqn:qr}), 
we see that all $F_{n,m}(\pi)(s_i)$'s
are noncommutative power series (=NPS's)
in $\pi(r_i)$'s and $\pi(r_i^*)$'s
for any $\pi\in\rep\co{m}$.
From (\ref{eqn:fminf}) and (\ref{eqn:finfncal}),
all $F_{\infty,n}(\pi)(t_j)$'s
are NPS's
in $\pi(s_i)$'s and $\pi(s_i^*)$'s
for any $\pi\in\rep\con$.
From (\ref{eqn:dashinf}) and (\ref{eqn:pidashfour}),
all $F_{n,\infty}(\pi)(s_i)$'s
are NPS's
in $\pi(t_j)$'s and $\pi(t_j^*)$'s
for any $\pi\in\rep\coni$.

\noindent
(vi)
From definitions of $F_{n,m}$'s,
the statement holds. 
\qedh

%
%
\begin{rem}
\label{rem:first}
{\rm
\begin{enumerate}
\item
From Theorem \ref{Thm:first}(iv),
we can prove that 
$F_{\infty,n}$ in (\ref{eqn:finfncal})
is essentially surjective, that is,
for any $\pi$ in $\rep\coni$,
there exists $\Pi$ in $\rep\con$ such that $F_{\infty,n}(\Pi)$ and $\pi$
are unitarily equivalent.
On the other hand,
for $\pi_1,\pi_2$ in Lemma \ref{lem:diff},
$F_{\infty,n}(\pi_1)$ and 
$F_{\infty,n}(\pi_2)$ are unitarily equivalent,
but not are $\pi_1$ and $\pi_2$.
Hence $F_{\infty,n}:\rep\con\to\rep\coni$ is not an isomorphism.
From this,
$F_{n,\infty}\circ F_{\infty,n}\ne id_{\rep\con}$
and 
$F_{n,\infty}\circ F_{\infty,m}\ne F_{n,m}$.
\item
Since $\con$ and $\co{m}$ are not strongly Morita equivalent
when $n\ne m$,
there exists no imprimitivity bimodule for $\con$ and $\co{m}$
corresponding to $F_{n,m}$.
\end{enumerate}
}
\end{rem}

%
%
\sftt{Formulas of $F_{n,m}$}
\label{section:fifth}
In order to check Theorem \ref{Thm:first},
we show concrete examples.
%
%
%
\ssft{General formulas}
\label{subsection:fifthone}
First, we show general formulas of $F_{n,m}$
in Theorem \ref{Thm:first}(v).
For $2\leq n,m<\infty$,
let $s_1,\ldots,s_n$ and $r_1,\ldots,r_m$
denote Cuntz generators of $\con$ and $\co{m}$,
respectively.
In general,
for any $\pi\in\rep\co{m}$,
we obtain
%
%
\begin{equation}
\label{eqn:rmiq}
\pi(r_{m})(I-Q_{\pi})
=
\pi(r_m)-
\sum_{l=0}^{\infty}
\sum_{i=1}^{m-1}\pi(r_m^{l+1}\,r_ir_i^*\,(r_{m}^l)^*)
\end{equation}
from the definition of $f_{m,\infty}$.
%
%
\begin{prop}
\label{prop:formal}
Fix $2\leq n,m<\infty$.
For $\pi\in \rep\co{m}$, 
let $\Pi:=F_{n,m}(\pi)\in\rep\con$.
Then the following hold.
\begin{enumerate}
\item
If $n<m$, then
%
%
\begin{equation}
\label{eqn:fmnasjleft}
\Pi(s_i)=
\left\{
\begin{array}{ll}
\pi(r_i) &(i<n),\\
\\
\disp{
\pi(r_m)(I-Q_{\pi})
+
\sum_{l=0}^{\infty}
\pi(\, r_m^l\,A_{n,m}\,(r_m^l)^*\, )}
&(i=n)
\end{array}
\right.
\end{equation}
where $A_{n,m}\in \co{m}$ is defined as
%
%
\begin{equation}
\label{eqn:anm}
A_{n,m}:=
\sum_{i=1}^{m-n}
r_{i+n-1}\,r_i^*
+
\sum_{i=m-n+1}^{m-1}
r_mr_{i+n-m}\,r_i^*.
\end{equation}
%
\item
Assume $n\geq m$.
\begin{enumerate}
%
\item
If $n\equiv 1 \mod{m-1}$, then 
%
%
\begin{equation}
\label{eqn:pisjand}
\begin{array}{rl}
\Pi(s_{j})=&
\left\{
\begin{array}{ll}
\pi(r_m^lr_i)\quad& 
\left(j< n
\right),\\
\\
\pi(r_m)(I-Q_{\pi})+\pi(r_m)^{k_0}Q_{\pi}\quad&(j=n)
\end{array}
\right.
\end{array}
\end{equation}
where 
$l$ and $i$ in the case of $j< n$
are defined as
$j=(m-1)l+i$ and $1\leq i\leq m-1$,
and
$k_0\geq 1$ in the case of $j= n$
is defined as $n=(m-1)k_0+1$.
\item
If $n\not \equiv 1 \mod{m-1}$, then 
%
%
\begin{equation}
\label{eqn:fanmipsjbeign}
\Pi(s_{j})=
\left\{
\begin{array}{ll}
\pi(r_m^lr_i)\quad &
\left(j< n
\right),
\\
\\
\disp{\pi(r_m)(I-Q_{\pi})+
\sum_{k=0}^{\infty}
\pi(\, r_m^{k+k_0}A_{n,m}\, (r_m^k)^*\, )}
\quad&(j=n)
\end{array}
\right.
\end{equation}
where $l$ and $i$ 
in the case of $j< n$
are defined as
$j=(m-1)l+i$ and
$1\leq i\leq m-1$,
and
$A_{n,m}\in \co{m}$ is defined as
%
%
\begin{equation}
\label{eqn:amnsum}
A_{n,m}:=
\sum_{i=1}^{m-j_0}
r_{i+j_0-1}r_i^*
+
\sum_{i=m-j_0+1}^{m-1}
r_mr_{i+j_0-m}r_i^*,
\end{equation}
and 
$k_0\geq 1$
and $2\leq j_0\leq m-1$
in the case of $j= n$
are defined as
$n=(m-1)k_0+j_0$.
\end{enumerate}
\end{enumerate}
\end{prop}
%
%
\pr
Let $t_1,t_2,\ldots$ denote Cuntz generators of $\coni$.
We write $\pi(r_i)$ as $r_i$ for short.
Since $\Pi(s_n)=r_m\,(I-Q_{\pi})+R_{\pi,n-1}$,	
the nontrivial part is the formula of $R_{\pi,n-1}$.
By definition,
%
%
\begin{equation}
\label{eqn:PIeqn}
\begin{array}{rl}
R_{\pi,n-1}
=&\disp{
\sum_{j=1}^{\infty}
f_{m,\infty}(t_{j+n-1}t_j^*)
}\\
\\
=&\disp{
\sum_{l=0}^{\infty}
\sum_{i=1}^{m-1}
f_{m,\infty}(t_{(m-1)l+i+n-1})\,f_{m,\infty}(t_{(m-1)l+i})^*.
}\\
\end{array}
\end{equation}

\noindent
(i)
For $i=1,\ldots,n-1$,
$\Pi(s_i)=f_{m,\infty}(t_i)=r_i$
because $1\leq i\leq n-1<n< m$.
On the other hand,
we see $f_{m,\infty}(t_{(m-1)l+i})=r^l_mr_i$
and
%
%
\begin{equation}
\label{eqn:ftmminus}
\begin{array}{rl}
f_{m,\infty}(t_{(m-1)l+i+n-1})
=&\disp{
\left\{
\begin{array}{ll}
r_m^lr_{i+n-1}\quad&(i+n-1\leq m-1),\\
\\
r_m^{l+1}r_{i+n-1-(m-1)}\quad&(i+n-1\geq m)\\
\end{array}
\right.}\\
\\
=&\disp{
\left\{
\begin{array}{ll}
r_m^lr_{i+n-1}\quad&(i\leq m-n),\\
\\
r_m^{l+1}r_{i+n-m}\quad&(i\geq m-n+1).\\
\end{array}
\right.}\\
\\
\end{array}
\end{equation}
Therefore
%
%
\begin{equation}
\label{eqn:mboxref}
\begin{array}{rl}
\mbox{(\ref{eqn:PIeqn})}
=&
\disp{
\sum_{l=0}^{\infty}
\left\{
\sum_{i=1}^{m-n}
r_m^lr_{i+n-1}
r_i^*(r_m^l)^*
+
\sum_{i=m-n+1}^{m-1}
r_m^{l+1}r_{i+n-m}
r_i^*(r_m^l)^*
\right\}
}\\
\\
=&
\disp{
\sum_{l=0}^{\infty}
r_m^l
\left\{
\sum_{i=1}^{m-n}
r_{i+n-1}
r_i^*
+
\sum_{i=m-n+1}^{m-1}
r_mr_{i+n-m}
r_i^*
\right\}
(r_m^l)^*.
}
\end{array}
\end{equation}
Hence the statement holds.

\noindent
(ii)
We see $\Pi(s_{(m-1)k+i})=f_{m,\infty}(t_{(m-1)k+i})=r_m^kr_i$
when $1\leq (m-1)k+i\leq n-1$.
Hence the case of $j< n$ is proved.
We show the case of $\Pi(s_n)$ as follows:

\noindent
(ii-a)
Since $n=(m-1)k_0+1$,
%
%
%
\begin{equation}
\label{eqn:fninftnumkia}
\begin{array}{rl}
f_{m,\infty}(t_{(m-1)k+i+n-1})
=&f_{m,\infty}(t_{(m-1)k+i+\{ (m-1)k_0+1\}-1})\\
=&f_{m,\infty}(t_{(m-1)(k+k_0)+i})\\
=&r_m^{k+k_0}r_i.
\end{array}
\end{equation}
From this, we obtain
%
%
\begin{equation}
\label{eqn:mboxrefsb}
\begin{array}{rl}
\mbox{(\ref{eqn:PIeqn})}
=&\disp{
\sum_{k=0}^{\infty}
\sum_{i=1}^{m-1}
r_m^{k+k_0}r_{i}r_i^*(r_m^k)^*
}
=r_m^{k_0}Q_{\pi}.
\end{array}
\end{equation}
Hence the statement holds.

\noindent
(ii-b)
Since $n=(m-1)k_0+j_0$,
%
%
\begin{equation}
\label{eqn:fninftnumkic}
\begin{array}{rl}
f_{m,\infty}(t_{(m-1)k+i+n-1})
=&f_{m,\infty}(t_{(m-1)k+i+\{ (m-1)k_0+j_0\}-1})\\
=&f_{m,\infty}(t_{(m-1)(k+k_0)+i+j_0-1})\\
\\
=&\disp{
\left\{
\begin{array}{ll}
r_m^{k+k_0}r_{i+j_0-1}\quad & (i+j_0-1\leq m-1),\\
\\
r_m^{k+k_0+1}r_{i+j_0-1-(m-1)}\quad & (i+j_0-1\geq  m)\\
\end{array}
\right.}\\
\\
=&\disp{
\left\{
\begin{array}{ll}
r_m^{k+k_0}r_{i+j_0-1}\quad & (i\leq m-j_0),\\
\\
r_m^{k+k_0+1}r_{i+j_0-m}\quad & (i\geq  m+1-j_0).\\
\end{array}
\right.}\\
\\
\end{array}
\end{equation}
\noindent
From this, we obtain
%
%
\begin{equation}
\label{eqn:mboxrefs}
\begin{array}{rl}
\mbox{(\ref{eqn:PIeqn})}
=&\disp{
\sum_{k=0}^{\infty}
\left\{
\sum_{i=1}^{m-j_0}
r_m^{k+k_0}r_{i+j_0-1}r_i^*(r_m^k)^*
+
\sum_{i=m-j_0+1}^{m-1}
r_m^{k+k_0+1}r_{i+j_0-m}r_i^*(r_m^k)^*
\right\}
}\\
=&\disp{
\sum_{k=0}^{\infty}
r_m^{k+k_0}
\left\{
\sum_{i=1}^{m-j_0}
r_{i+j_0-1}r_i^*
+
\sum_{i=m-j_0+1}^{m-1}
r_mr_{i+j_0-m}r_i^*
\right\}(r_m^k)^*.
}
\end{array}
\end{equation}
Hence the statement holds.
\qedh

\noindent
In Proposition \ref{prop:formal}(ii-a),
if $n=m$, that is, $k_0=1$, then $F_{n,n}(\pi)=\pi$ is reproved.

%
%
\ssft{$F_{2,3}$ and $F_{3,2}$}
\label{subsection:fifthtwo}
We show the case of $(n,m)=(2,3)$
in Proposition \ref{prop:formal}(i)
and
the case of $(n,m)=(3,2)$, that is, $3=(2-1)\cdot 2+1$,
in Proposition \ref{prop:formal}(ii-a) with $k_0=2$.
Let 
$s_1,s_2$ and  $r_1,r_2,r_3$ denote Cuntz generators of $\co{2}$ and
$\co{3}$, respectively.

For $\pi\in\rep\co{3}$,
we obtain $F_{2,3}:\rep\co{3}\to \rep\co{2}$ as 
%
%
\begin{equation}
\label{eqn:ftwothree}
\left\{
\begin{array}{rl}
F_{2,3}(\pi)(s_1)=&\pi(r_1),\\
\\
F_{2,3}(\pi)(s_2)
=&
\disp{\pi(r_3)
-\sum_{k=0}^{\infty}
\pi(r_3^{k+1}(r_1r_1^*+r_2r_2^*)(r_3^k)^*)
}\\
&+\disp{
\sum_{k=0}^{\infty}
\pi(r_3^k(r_2r_1^*+r_3r_1r_2^*)(r_3^k)^*)}\\
\\
=&
\disp{\pi(r_3)+
\sum_{k=0}^{\infty}
\pi(\,
r_3^k\,\{r_2r_1^*+r_3(r_1r_2^*-r_1r_1^*-r_2r_2^*)\}\,(r_3^k)^*\, ).}
\end{array}
\right.
\end{equation}

\noindent
For $\pi\in\rep\co{2}$,
we obtain $F_{3,2}:\rep\co{2}\to \rep\co{3}$ as 
%
%
\begin{equation}
\label{eqn:ftwothreeb}
\left\{
\begin{array}{rl}
F_{3,2}(\pi)(r_1)=&\pi(s_1),\\
\\
F_{3,2}(\pi)(r_2)=&\pi(s_2s_1),\\
\\
F_{3,2}(\pi)(r_3)
=&
\disp{\pi(s_2)
-\sum_{k=0}^{\infty}
\pi(s_2^{k+1}s_1s_1^*(s_2^k)^*)
}+\disp{\sum_{k=0}^{\infty}
\pi(s_2^{k+2}s_1s_1^*(s_2^k)^*)}\\
\\
=&
\disp{\pi(s_2)
+\sum_{k=0}^{\infty}
\pi(s_2^{k+1}(s_2-I)s_1s_1^*(s_2^k)^*).}
\end{array}
\right.
\end{equation}
In addition,
if $\pi$ satisfies $Q_{\pi}=I$, then 
$F_{3,2}(\pi)(r_3)=\pi(s_2^2)$.
Furthermore,
for any $2\leq n,m<\infty$,
if $\pi$ satisfies $Q_{\pi}=I$, then
we see that
all $F_{n,m}(\pi)(s_i)$'s in Proposition \ref{prop:formal}(ii-a),
that is, $n=(m-1)k_0+1$ for some $k_0\geq 1$,
are written as noncommutative polynomials in
$\pi(r_i)$'s and $\pi(r_i^*)$'s.
This implies
that there exists $f_{m,n}\in \hom(\con,\co{m})$
such that $F_{n,m}(\pi)=\pi\circ f_{m,n}$
for any $\pi\in\rep\co{m}$.
%
%
%
\ssft{$F_{3,4}$ and $F_{4,3}$}
\label{subsection:fifththree}
We show
the case of $(n,m)=(3,4)$
in Proposition \ref{prop:formal}(i)
and
the case of $(n,m)=(4,3)$, that is, $4=(3-1)\cdot 2+2$,
in Proposition \ref{prop:formal}(ii-b) with $(k_0,j_0)=(2,2)$.
Let $s_1,s_2,s_3,s_4$ and $r_1,r_2,r_3$ denote
Cuntz generators of $\co{4}$ and $\co{3}$, respectively.

For $\pi\in\rep\co{4}$,
we obtain $F_{3,4}:\rep\co{4}\to \rep\co{3}$ as
%
%
\begin{equation}
\label{eqn:leftftrhree}
\left\{
\begin{array}{rl}
F_{3,4}(\pi)(r_i)=&\pi(s_i)\quad(i=1,2),
\\
\\
F_{3,4}(\pi)(r_3)
=&
\disp{
\pi(s_4)
-\sum_{k=0}^{\infty}
\pi(\, s_4^{k+1}(s_1s_1^*+s_2s_2^*+s_3s_3^*)(s_4^k)^*\, )
}\\
&+\disp{\sum_{k=0}^{\infty}
\pi(\, s_4^k(s_{3}s_{1}^*
+s_4s_{1}s_{2}^*
+s_4s_{2}s_{3}^*)(s_{4}^k)^*\, )}\\
\\
=&
\disp{
\pi(s_4)}\\
&\disp{
+\sum_{k=0}^{\infty}
\pi(\, s_4^k\{s_{3}s_{1}^*
+s_4(s_{1}s_{2}^*+s_{2}s_{3}^*-I+s_4s_4^*)\}(s_{4}^k)^*\, )}
\end{array}
\right.
\end{equation}
because 
$s_1s_1^*+s_2s_2^*+s_3s_3^*=I-s_4s_4^*$.
For $\pi\in\rep\co{3}$,
we obtain $F_{4,3}:\rep\co{3}\to \rep\co{4}$ as
%
%
\begin{equation}
\label{eqn:ffourthree}
\left\{
\begin{array}{rl}
F_{4,3}(\pi)(s_i)=& \pi(r_i)\quad(i=1,2),\\
\\
F_{4,3}(\pi)(s_3)=& \pi(r_3r_1),\\
\\
F_{4,3}(\pi)(s_4)
=&\disp{
\pi(r_3)-
\sum_{k=0}^{\infty}
\pi(\, r_3^{k+1}(r_1r_1^*+r_2r_2^*)(r_3^k)^*\, )
}\\
&+\disp{
\sum_{k=0}^{\infty}
\pi(\, r_3^{k+2}
(r_2r_1^*+r_3r_1r_2^*)(r_3^k)^*\, )}\\
\\
=&
\disp{
\pi(r_3)+\sum_{k=0}^{\infty}
\pi(\, r_3^{k+1}
(r_3r_2r_1^*+r_3^2r_1r_2^*-r_1r_1^*-r_2r_2^*)(r_3^k)^*\, ).}
\end{array}
\right.
\end{equation}

\noindent
Remark 
$\hom(\co{2},\co{3})=
\hom(\co{3},\co{4})=\hom(\co{4},\co{3})=\emptyset$.
In general, the following holds.
%
%
\begin{prop}
\label{prop:ktheory}
{\rm (\cite[Lemma 2.1]{SE01})}
For $2\leq m,n<\infty$, ${\rm Hom}(\co{m},\con)\ne \emptyset$ if and only if
there exists a positive integer $k$ such that $m=(n-1)k+1$, that is,
$m\geq n$ and $m\equiv 1\mod n-1$,
where
$\hom(A,B)$ denotes the 
set of all unital $^*$ homomorphisms from $A$ to $B$.
\end{prop}

These types of power series presentations 
have been obtained for other algebras 
in \cite{IWF01,RBS01,BFO01,UFB01}.
Concrete power series presentations are useful
to construct interesting operators like $R_{\pi,a}$ in (\ref{eqn:qr}).
For example, see \cite[(1.6)]{RBS01}.

%
%
\sftt{Open problems}
\label{section:sixth}
\begin{enumerate}
\item
Find an example of two associative unital nontopological rings such that
they are not isomorphic but Morita isomorphic.
For example,
for a ring  $A$ and $n\geq 2$,
are $A$ and $M_n(A)$ Morita isomorphic?
A one-sided invertible functor can be obtained by standard way,
but the author fails to construct a two-sided invertible functor
for their module categories.
It seems that 
an algebraic version of Theorem \ref{Thm:first}(vi) is necessary
in order to find such an example.
\item
For any $2\leq n<\infty$,
it is known that $\coni$ and $\con$ are Morita equivalent.
Are $\coni$ and $\con$ Morita isomorphic?
(See Remark \ref{rem:first}(i).)
\item
Find an example of 
two non-type I, separable, nuclear $C^*$-algebras 
$A$ and $B$
such that $A$ and $B$ are not Morita isomorphic.
(See Theorem \ref{Thm:Beer}.)
\item
For a pair $(A,B)$ of $C^*$-algebras,
find a necessary and sufficient condition 
that $A$ and $B$ are Morita isomorphic.
For example,
for a pair $(A,B)$ of Cuntz-Krieger algebras (\cite{CK}),
which are well-known generalizations of Cuntz algebras,
find a necessary and sufficient condition 
that $A$ and $B$ are Morita isomorphic.
About related topics, see \cite{Matsumoto}.
\item
When $n\ne m$,
$\con$ and $\co{m}$ are 
Morita isomorphic but not stably isomorphic.
Show a relation between Morita isomorphism and stably isomorphism
for general $C^*$-algebras.
\item
A main philosophy of mathematics 
is an identification up to isomorphism.
Therefore category equivalence (=CE) just represents this idea faithfully.
According to \cite[Chap IV, $\S$4]{Maclane},
the notion of CE 
is more general and more useful than 
that of category isomorphism (=CI).
Duality theorems in functional analysis are often
instances of CE's.
In comparison to CE, studies of CI are very few (\cite{Traylor}).
Under like this unfavorable aspect,
can we find a new significance of CI?
\end{enumerate}

\noindent
{\bf Acknowledgments.}
The author would like to express his sincere thanks to 
Professor Satoshi Yamanaka
for discussion about Morita equivalence at Ajinomingei.

\appendix
%
%
\sftt{A short survey of Morita equivalence}
\label{section:appone}
We outline Morita equivalences for rings and $C^*$-algebras.
%
%
\ssft{Morita equivalences for rings}
\label{subsection:apponeone}
Two associative unital rings are said to be {\it Morita equivalent}
if categories of their left modules are equivalent (\cite{Morita1}).
This is equivalent to an existence of progenerator 
for their module categories
by  {\it Morita's basic theorem} (\cite{Morita1}, 
see also \cite[(18.26) Theorem]{Lam}).
Hence a Morita equivalence 
between two rings implies 
an existence of correspondence of their modules, 
which is interesting as a method of construction of module.
For example,
a ring $A$ and the full matrix ring $M_n(A)$ with components in $A$
are Morita equivalent by the functor $A$-${\rm Mod}\ni M\mapsto 
F(M):=M^{\oplus n}
\in M_n(A)$-${\rm Mod}$.
Then $F$ is isomorphic to
the functor
$P\otimes_A-$ for the progenerator $P:=A^{\oplus n}$.
For two commutative rings  $A$ and $B$,
they are Morita equivalent if and only if 
they are isomorphic 
(\cite[p522]{Aluffi}).
Therefore 
$A$ and $B$ are Morita equivalent if and only if they
are Morita isomorphic.
About applications of Morita equivalence, see \cite{AF,Yamanaka}.
%
%
\ssft{Two Morita equivalences for $C^*$-algebras}
\label{subsection:apponetwo}
On the other hand,
for topological rings,
it is more complicated because topologies are necessary for their modules.  
By only the straightforward generalization,
Morita's basic theorem does not hold for $C^*$-algebras.
Therefore it is necessary to distinguish an equivalence of module categories
and an existence of C$^*$-version of progenerator.

For $C^*$-algebras,
Rieffel introduced two notions
about Morita equivalence:
Two $C^*$-algebras are said to be 
{\it Morita equivalent} (\cite[8.17]{Rieffel1974a})
if their categories of nondegenerate 
$^*$ representations are equivalent;
Two $C^*$-algebras 
are said to be 
{\it strongly Morita equivalent}
(\cite[Proposition 6.26]{Rieffel1974b})
if there exists an imprimitivity 
bimodule for them
(see also \cite[13.7.1]{BlackadarK} and  \cite[Chapter 7]{Lance}),
which is a $C^*$-version of progenerator.
Remark that the definition of Morita equivalence 
for rings and that for $C^*$-algebra are different.
Furthermore,
Morita equivalence for $C^*$-algebras is often 
used as a meaning of strong Morita equivalence (\cite{Blackadar}).
We review known results as follows.
%
%
\begin{Thm}
\label{Thm:Beer}
\begin{enumerate}
\item
{\rm (\cite[3.7. Theorem]{Beer})}
All non-type {\rm I}, separable, nuclear $C^*$-algebras are 
Morita equivalent.
\item
{\rm (\cite{BGR}, see also \cite[II.6.6.12]{Blackadar})}
If two $\sigma$-unital $C^*$-algebras 
are strongly Morita equivalent, then they are stably isomorphic.
\end{enumerate}
\end{Thm}
From Theorem \ref{Thm:Beer}(i),
all Cuntz algebras are Morita equivalent.
From Theorem \ref{Thm:Beer}(ii), 
two Cuntz algebras are not strongly Morita equivalent
when they are not isomorphic.
Two unital $C^*$-algebras are strongly
Morita equivalent if and only if
they are Morita equivalent as rings (\cite[$\S$1.8, Theorem]{Beer}).
About basic properties of strong Morita equivalence, see \cite{Blackadar}.
About related topics, see \cite{BM,BMP,RW}.

%
%

\label{Lastpage}

\end{document}